\documentclass[12pt]{article}
\usepackage{graphicx}
\usepackage{fullpage}
\begin{document}

\section*{Proof Without Words}

\begin{minipage}[t]{3in}
Japheth Wood\\
Bard College\\
Annandale-on-Hudson, NY 12504
\end{minipage}
\begin{minipage}[c]{3in} \centering
\includegraphics[scale=1.5,trim=31 24 35 32,clip]{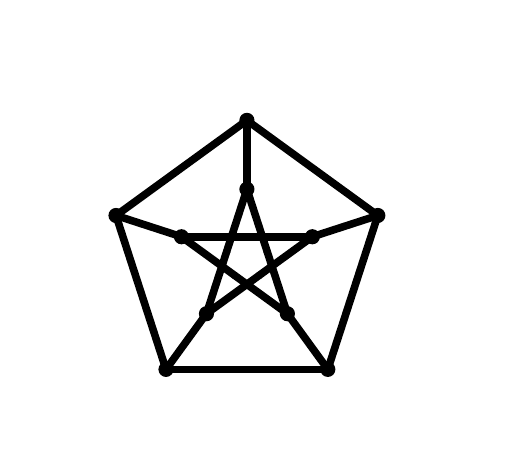}

The Petersen Graph
\end{minipage}

\begin{center}
	\includegraphics[scale=2.3,trim=30 26 74 26,clip]{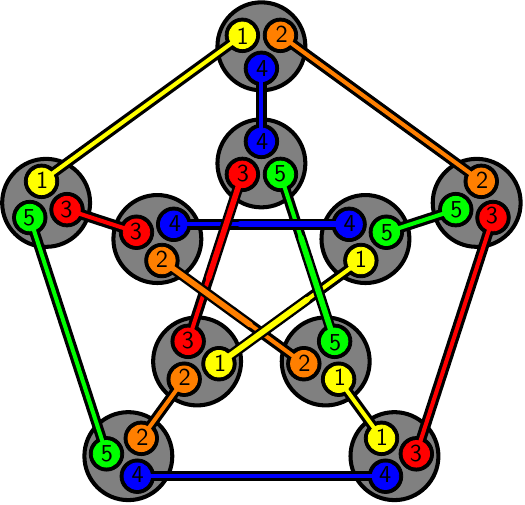}
	
	The Automorphism Group of the Petersen Graph is isomorphic to $S_5$.
\end{center}


\subsection*{Abstract}
The automorphism group of the Petersen Graph is shown to be isomorphic to the symmetric group on 5 elements. The image represents the Petersen Graph with the ten 3-element subsets of $\{1, 2, 3, 4, 5\}$ as vertices. Two vertices are adjacent when they have precisely one element in common.
\end{document}